\documentclass{IEEEtran}
\usepackage{amsmath,amsfonts}
\usepackage{algorithmic}
\usepackage{algorithm}
\usepackage{array}
\usepackage[caption=false,font=normalsize,labelfont=sf,textfont=sf]{subfig}
\usepackage{textcomp}
\usepackage{stfloats}
\usepackage{url}
\usepackage{verbatim}
\usepackage{graphicx}
\usepackage{cite}
\usepackage{amssymb}
\usepackage{amsthm}
\usepackage{color}
\usepackage{arydshln}

\begin{document}

\title{Data-driven $H_{\infty}$ predictive control for constrained systems: a Lagrange duality approach}

\author{Wenhuang Wu,~
        Lulu Guo,~   
        Nan Li,~\IEEEmembership{Member,~IEEE,}   
        Hong Chen,~\IEEEmembership{Fellow,~IEEE}   
\thanks{This research was supported in part by the National Natural Science Foundation of China under Grant No. 62203334, 
Shanghai Rising-Star Program under Grant No. 22YF1451300, and the Fundamental Research Funds for the Central Universities. \emph{(Corresponding author: Lulu Guo.)}}
\thanks{Wenhuang Wu is with the Shanghai Research Institute for Intelligent Autonomous Systems, Tongji University, Shanghai 200092, China (e-mail: wuwh898@tongji.edu.cn).}
\thanks{Lulu Guo and Hong Chen are with the Department of Control Science and Engineering, Tongji University, Shanghai 201804, China (e-mail: guoll21@tongji.edu.cn; chenhong2019@tongji.edu.cn).}
\thanks{Nan Li is with the School of Automotive Studies, Tongji University, Shanghai 201804, China (e-mail: li\_nan@tongji.edu.cn).}}

\maketitle

\begin{abstract}
  This article proposes a data-driven $H_{\infty}$ control scheme for time-domain constrained systems based on model predictive control formulation.
  The scheme combines $H_{\infty}$ control and minimax model predictive control, enabling a more effective way to handle external disturbances and time-domain constraints.
  First, by leveraging input-output-disturbance data, the scheme ensures $H_{\infty}$ performance of the closed-loop system. 
  Then, a minimax optimization problem is converted into a more manageable minimization problem employing Lagrange duality, 
  which reduces conservatism typically associated with ellipsoidal evaluations of time-domain constraints.
  The study examines key closed-loop properties, including stability, disturbance attenuation,
  and constraint satisfaction, achieved by the proposed data-driven moving horizon predictive control algorithm. 
  The effectiveness and advantages of the proposed method are demonstrated through numerical simulations involving a batch reactor system, 
  confirming its robustness and feasibility under noisy conditions.
\end{abstract}


\section{Introduction}
\IEEEPARstart{W}{ith} the development of data science and artificial intelligence, the analysis and control of systems based on data-driven control frameworks have become increasingly popular 
in recent years. 
This trend is attributable to the abundance of easily collected system data, which enables bypassing the challenges associated with creating highly accurate models.
Different from traditional control methods that rely on predefined system models and rule-based algorithms, 
which may struggle to demonstrate good control performance due to the inherent variability in certain complex scenarios, 
data-driven control is derived directly from the data generated by the system itself \cite{de2019formulas,seuret2024robust,berberich2020robust,van2020noisy,10415858}. 
Control strategies based on data-driven control are widely applied in many fields, including power systems, biomedical engineering, and robotic control \cite{zhang2015data,sun2021machine,novara2019data,bruder2020data}. 
Extensive research has been conducted on the analysis and control of systems using data-driven control frameworks, 
including robust control \cite{karimi2017data}, nonlinear control \cite{tanaskovic2017data}, 
optimal control \cite{martinelli2022data}, and model predictive control (MPC) \cite{berberich2020data}.

MPC has emerged as a leading control strategy in the realm of advanced process control due to its robust ability to handle multivariable systems and constraints. 
At its essence, MPC is an optimization-based control approach that utilizes a dynamic model of the system to predict and optimize future behavior over a finite time horizon\cite{kouvaritakis2016model}.
Hence, the operation of MPC usually depends on a system model, which can be derived from first principles or system identification using available data \cite{ljung1987theory}. 
Recently, there has been a rise in the popularity of MPC controller design methods based on data-driven control, which utilize the available data to solve 
an optimization problem that determines the optimal control actions \cite{ABOLPOUR2024225,carron2019data,jafarzadeh2021dmpc,smarra2018data,han2019knowledge,bongard2022robust,piga2017direct,rosolia2017learning}.
Typical examples of these methods include learning-based or adaptive MPC \cite{HARA2023111271,10374442}, MPC based on Koopman operator \cite{KORDA2018149}, 
MPC based on behavioral systems theory \cite{9303965}, MPC based on noisy data \cite{xie2024datadriven}.
In \cite{liu2022predictor}, a novel adaptive data-driven MPC approach for power converters has been proposed, integrating neural network-based predictors and finite control set MPC. 
To mitigate uncertainties, a supervised imitation learning technique transforms most of the online computational burden to offline computation, 
using a trained artificial neural network to enhance robustness and simplify implementation in practical applications.
In \cite{zhang2022robust}, a robust MPC approach for nonlinear discrete-time dynamical systems using Koopman operators has been presented. 
The proposed method combines a nominal MPC using a lifted Koopman model with an offline nonlinear feedback policy to ensure 
closed-loop robustness against modeling errors and external disturbances, while also guaranteeing convergence properties.
In \cite{pan2022stochastic}, the extension of the deterministic fundamental lemma to stochastic systems using polynomial chaos expansions has been discussed, 
This method allows the prediction of future probability distributions for a linear time-invariant (LTI) system with random parameters based on previously recorded data and disturbance distributions. 
Building on this extension, a conceptual framework for data-driven stochastic predictive control has been introduced.

Moreover, real-world systems often encounter uncertainties, disturbances, and model inaccuracies that can compromise the performance of standard MPC.
To address these challenges, Robust Model Predictive Control (RMPC) has been developed.
RMPC extends the traditional MPC framework by incorporating robustness into the optimization process, 
ensuring reliable performance even under uncertain conditions. By considering worst-case scenarios and employing robust optimization techniques, 
RMPC maintains system stability and performance despite external disturbances and parameter variations \cite{bemporad2007robust,lofberg2003minimax}.
To cope with situations where the system model cannot be accurately obtained, RMPC based on data-driven control is a feasible approach. 
Currently, research on data-driven RMPC (DDRMPC) has been conducted from various perspectives. For example, the stability and robustness guarantees 
based on DDRMPC are studied in \cite{berberich2020robust}, and DDRMPC utilizing Koopman operators is presented in \cite{mamakoukas2022robust}.
Furthermore, there are other recent studies on DDRMPC methods \cite{chen2021semiclosed,mahmood2023data,bongard2022robust}. 
The more specific results are as follows:
in \cite{li2024data}, a robust data-driven predictive control method for constrained systems has been presented, 
integrating MPC and $H_{\infty}$ control for addressing constraints and disturbance rejection. 
This method adjusts $H_{\infty}$ performance online by using system states and disturbance forecasts, ensuring robust stability, 
constraint satisfaction, and disturbance attenuation even subject to disturbance conditions.
In \cite{zhang2024data}, an event-triggered robust MPC design method for unknown systems using initially measured input-output data has been proposed. 
The method introduces terminal inequality constraints to enhance the feasible region and reduce the prediction horizon, 
ensuring recursive feasibility and input-to-state stability under mild conditions, with simulation results demonstrating its effectiveness.
In \cite{deng2024event}, a novel data-driven robust iterative learning predictive control scheme for multiple-input multiple-output nonaffine nonlinear systems with actuator constraints has been presented. 
The proposed approach leverages a noise-tolerant zeroing neural network for adaptive estimation and demonstrates effective multivariable tracking and noise suppression,
which is validated through a proton exchange membrane fuel cell thermal management system.

In this context, 
this paper not only addresses stability but also leverages $H_{\infty}$ control and MPC to tackle disturbance suppression and constraint satisfaction 
in data-driven control. Regarding data-driven $H_{\infty}$ control, \cite{van2020noisy} introduces a method based on noise-corrupted data, 
which effectively addresses the common issues of high dimensionality and robustness degradation caused by noise in data-driven control. 
However, the approach in \cite{van2020noisy} neither considers time-domain constraints nor achieves dynamic optimization performance. 
These limitations are addressed in \cite{li2024data}, where a data-driven $H_{\infty}$ MPC method is proposed. 
Nevertheless, when the system encounters more disturbances—such as disturbances in the output or when the system disturbance matrix 
is not an identity matrix—the control method in \cite{li2024data} fails to perform effectively. 
To address this gap, this article investigates an $H_{\infty}$ predictive control scheme based on a data-driven framework, 
providing a more comprehensive solution to disturbance suppression and constraint handling in data-driven control.
The scheme employs a Lagrange duality approach to convert a minimax optimization into a minimization problem, 
which alleviates the challenges induced by infeasibility and time-domain constraints. 
By utilizing input-state and disturbance data, the closed-loop guarantees of controlled systems based on a moving horizon predictive 
control algorithm, including the satisfaction of time-domain constraints, stability, and $H_{\infty}$ control performance, can be achieved. 
The detailed contributions are summarized as follows:
 
(1) In this paper, we propose a novel $H_{\infty}$ predictive control scheme for constrained systems based on input-state-disturbance data. 
To address more general system scenarios, the proposed scheme leverages the matrix Finsler's lemma \cite{van2021matrix},
which offers less conservatism compared to the traditional S-Lemma \cite{li2024data},
to derive the criterion for stabilization and $H_{\infty}$ disturbance attenuation control of an unknown system subject to output constraints.

(2) The proposed scheme employs a Lagrange duality approach within a data-driven control framework to transform a minimax optimization problem into 
a more tractable minimization problem. This transformation not only reduces the conservatism introduced by ellipsoid constraint analysis 
but also facilitates feasibility analysis, making the overall control design more practical and efficient.

(3) The paper presents a moving horizon predictive control algorithm by introducing a dissipativity inequality constraint to the data-driven $H_{\infty}$ minimization problem, 
which can ensure sustained disturbance suppression over time. Then, the closed-loop properties are discussed. 
The properties involve recursive feasibility of optimal control problem,
stability of closed-loop system, adaptive disturbance attenuation performance of $H_{\infty}$ control, and constraint satisfaction of output.  

The organization of the paper is summarized as follows: 
In Section 2, we outline the issue addressed in this paper, including an assumption and a lemma.
Section 3 introduces our approach for synthesizing an $H_{\infty}$ controller of an unknown system utilizing offline data.
The second problem addressed in Section 3 involves the relaxation of the minimax problem using a Lagrange duality method, 
and then a tractable minimization problem based on a data-driven control framework is proposed.
In the third part of Section 3, we propose a data-driven moving horizon optimal control algorithm and present closed-loop properties of the control algorithm. 
Section 4 demonstrates the effectiveness and superiority of the control algorithm with a practical example. 
Finally, the paper is concluded in Section 5.

{\bf{Notations:}} The symbol $\mathbb{R}$ represents the real number space, and its superscripts represent the corresponding dimensions. 
For example, $\mathbb{R}^q$ and $\mathbb{R}^{p\times q}$ represent $q$-dimensional real vectors and $p\times q$ real matrices, respectively.
$||f(x)||$ denotes the Euclidean norm of the function $f(x)$. We use $I$ to signify the identity matrix of appropriate dimensions. 
Similarly, the $0$ in a matrix indicates a zero matrix of appropriate dimensions. 
For a matrix $N$, $N^{-1}$, $N^{+}$ and $N^{T}$ indicate its inverse, pseudo-inverse
and transpose, respectively. Additionally, $N>0$ ($N\geq 0$) implies that $N$ is positive (semi-)definite. 
The symbol $\ast$ stands for the symmetric term corresponding to a specific term in a symmetric matrix.

\section{Problem formulation}

Consider a discrete-time LTI system as follows:
\begin{subequations}\label{sys_1_sum}
\begin{align}{x}(t+1)&=A_{\nu }{x}(t)+B_{\nu } u(t)+E_{\nu } {w}(t),  \label{sys_1}\\
    y_1(t) &= C_1 x(t)+ D_1u(t) + E_1 w(t),   \label{out_1}\\
    y_2(t) &= C_2 x(t)+ D_2u(t)     \label{out_2}  \end{align}
\end{subequations}
where $x\in \mathbb{R}^{p}$ is the state, $u\in \mathbb{R}^{q}$ represents the input, and $w\in \mathbb{R}^{l} $ denotes an external disturbance. 
The matrices $A_{\nu } \in \mathbb{R}^{p\times p} $, $B_{\nu } \in \mathbb{R}^{p\times q} $ and $E_{\nu } \in \mathbb{R}^{p\times l} $
denote the state, input and disturbance matrices, 
respectively, which are unknown.; $C_1$, $C_2$, $D_1$, $D_2$ and $E_1$ are assumed to be known matrices. 
The $y_1\in \mathbb{R}^{q_1}$ is the performance output vector, 
and $y_2\in \mathbb{R}^{q_2}$ represents the output vector constrained by time-domain constraints as follows:
\begin{align}|y_{2m}(t)|\leq y_{2m,~\max},~~~\forall t\geq 0, ~~m=1,2,\cdots,q_2.  \label{z_constraint}  
\end{align}
The basic assumptions about the system are that $(A_{\nu }, B_{\nu })$ is stabilizable and $(C_1, A_{\nu })$ is observable.

In the article, the fundamental goal is to design a controller to achieve closed-loop internal stability for LTI system (\ref{sys_1}) 
based on a certain number of sampling sequence of states and inputs of the system. 
To do this, assuming 
\begin{equation}
    \begin{aligned} 
    \hat{X}&:=[x(0)\quad x(1)\quad\cdots\quad x(T)] , \\
    \hat{U}&:=[u(0)\quad u(1)\quad\cdots\quad u(T-1)] ,        \\
    \hat{W}&:=[w(0)\quad w(1)\quad\cdots\quad w(T-1)]  ,      \\
    \hat{X}_{-}&:=[x(0)\quad x(1)\quad\cdots\quad x(T-1)] ,   \\
    \hat{X}_{+}&:=[x(1)\quad x(2)\quad\cdots\quad x(T)]    
\end{aligned}
\label{data}  
\end{equation}
where $\hat{U}$, $\hat{W}$ and $\hat{X}_{-}$ denote the sampled data sequences of previous inputs, disturbances and states over time interval $t\in [0,T-1]$;
$\hat{X}_{+}$ stands for the sampled data sequence of next states over time interval $t\in [1,T]$. Hence, the following equation can be clearly derived:
\begin{align}\hat{X}_+=A_{\nu }\hat{X}_-+B_{\nu }\hat{U} +E_{\nu } \hat{W}. \label{sys_data} \end{align}
Note that the matrices $A_{\nu }$, $B_{\nu }$ are unknown, 
while $\hat{X}$, $\hat{U}$ and $\hat{W}$ are sampled. The following assumption on the disturbance $\hat{W}$ is introduced.

{\bf{Assumption 1.}} 
In this article, the disturbance samples $w(0), w(1),..., w(T-1)$, which are collected in the matrix $\hat{W}$, are assumed to satisfy
\begin{equation}
\sum_{t=0}^{T-1}||w(t)||^2\leqslant \alpha  , \label{assump_1}
\end{equation}
where $||w(t)||^2\leq \alpha_t$ for a given scalar $\alpha_t \geq 0$.
This implies that the disturbances are energy-bounded (with an overall bound $\alpha\geq 0$). As $T\rightarrow \infty$, we have $w\in L_2[0,~\infty)$. 

For a more general situation, we introduce a set $\Pi$ to denote all systems $(A,B,E)$ compatible with  
the measurement sequence $\hat{X}_-$, $\hat{X}_+$, $\hat{U}$ and $\hat{W}$; that is to say, we have a general system
\begin{align} \hat{X}_+=A\hat{X}_-+B\hat{U}+ E \hat{W}  \label{sys_data_2}
\end{align}
where $\hat{W}$ satisfies the Assumption 1. The set $\Pi$ is formulated as 
\begin{align} 
\Pi=\{(A,~B,~{E})|~(\ref{sys_data_2})\text{ holds}\}.
\end{align}
For the set $\Pi$, if we can design a controller to stabilize all systems within it, then the true system $(A_{\nu},~B_{\nu},~E_{\nu})$ can also be stabilized.

On the basis of (\ref{sys_data}) and (\ref{assump_1}), we can not only achieve the fundamental goal below (\ref{z_constraint}) 
but also extend the result to disturbance attenuation controller design with the help of the following lemma. 

{\bf{Lemma 1.}}\cite{van2021matrix} Consider symmetric matrices $G, H \in \mathbb{R}^{(p+k)\times (p+k)}$ partitioned as follows:
\begin{equation}G=\begin{bmatrix}G_{11}& G_{12} \\ G^T_{12} & G_{22}\end{bmatrix},\quad H=\begin{bmatrix}H_{11}&H_{12}\\
        H_{12}^T&H_{22}\end{bmatrix}     \notag  
\end{equation}
 where $G_{11},~H_{11} \in \mathbb{R}^{p\times p}$. Assume that \\
(1) $G_{12}=0$, $G_{22}\leq 0$;  \\
(2) $H_{22}\leq 0$ and $H_{11}-H_{12}H_{22}^{+}H_{12}^{T}=0;$  \\
(3) $\exists F$ such that $G_{11}+F^{T}G_{22}F>0$ and $H_{22}F=H^T_{12}$.
Then, we have that
\begin{equation*}\begin{bmatrix}I\\F\end{bmatrix}^T G\begin{bmatrix}I\\F\end{bmatrix}\geq0, ~~~\forall F\in\mathbb{R}^{p\times k}
\end{equation*} 
such that
\begin{equation*}\begin{bmatrix}I\\F\end{bmatrix}^T H\begin{bmatrix}I\\F\end{bmatrix}=0\end{equation*} 
if and only if there exists $\lambda \in \mathbb{R}$ such that $G-\lambda H\geq 0$. 

\textbf{Proof:} See Theorem 1 of \cite{van2021matrix}.
\vspace{\baselineskip}

\textbf{Remark 1:}
In data-driven control systems, traditional model-based controllability and observability analysis 
(which relies on explicit system matrices) is replaced by methods that directly utilize input-output data. 
Instead of constructing controllability and observability matrices using system dynamics, we can extract information from collected trajectories 
to infer these properties.
\vspace{\baselineskip}

The purpose of introducing the above lemma is to bridge the Lyapunov stability condition presented in the next section and measurement data 
(\ref{data}) to achieve the stability of the system (\ref{sys_1_sum}) and extend them to another controller design. In what follows, we 
turn our attention to the optimization problem with respect to the performance output $y_1(t)$, 
and design the optimal controller to achieve closed-loop stability and dissipation by combining the results discussed above. 

\section{Main results}
In this section, we first address the stabilization and $H_\infty$ control of unknown systems based on measurement data.
We then introduce the Lagrange dual formulation for a minimax optimization problem with respect to performance output. 
Subsequently, we formulate the optimization problem with performance and output constraints. 
Finally, we present a moving horizon predictive control algorithm and analyze the robust performance of the controlled system.  

\subsection{Data-driven stabilization and $H_{\infty}$ control} 
Through the above discussion, we know that the primary task is to find the appropriate controller gain $K$ that can stabilize any $(A,B,E)\in \Pi$.
For the stabilization and disturbance attenuation performance of the system, the $H_\infty$ performance level $\gamma$
from the external input $w(t)$ to the control output $y_1(t)$ shall be minimized. 
Based on Lyapunov stability theorem, define the Lyapunov function $V(t)=x^T(t) P x(t)$. Then, the inequality $x^T(t) P x(t) -(Ax(t)+Bu(t))^T P(Ax(t)+Bu(t))>0$ implies that
the stabilization of a closed-loop system can be achieved if there exist a matrix $P=P^T>0$ and a feedback gain $K$ such that
\begin{equation}
     P-A^T_U P A_U   > 0   \label{lya_stab}
    \end{equation}
where $A_U=A+BK$. According to Proposition 3.12 in \cite{scherer1999lecture}, for zero initial conditions, 
the system (\ref{sys_1_sum}) has $H_\infty$ performance
level $\gamma$ from $w(t)$ to $y_1(t)$ if and only if 
$\sum_{t=0}^{\infty}\gamma^2 w^T(t) w(t) \geq \sum_{t=0}^{\infty} y_1^T(t) y_1(t)$, 
which can be presented in linear matrix inequality (LMI) form as follows
\begin{equation}
    \left[\begin{array}{cc} P-A^T_UPA_U-C^T_U C_U & -A^T_UP E -C^T_U E_1 \\ - E^T PA_U-E^T_1 C_U & -E^T PE -E^T_1 E_1+\gamma^2 I 
    \end{array}\right]
            > 0   \label{H_infty_2}
\end{equation}
where $C_U=C_1+D_1 K$. From (\ref{H_infty_2}), it follows that
\begin{align}
    \begin{bmatrix} P-C^T_UC_U & -C^T_U E_1\\ -E^T_1 C_U & \gamma^2I-E^T_1 E_1 \end{bmatrix}
    - \begin{bmatrix}  A^T_U P A_U&A^T_U P E \\E^T P A_U &  E^T P E \end{bmatrix}  > 0. \label{H_infty_3}
\end{align}
Furthermore, (\ref{H_infty_3}) is equivalent to
\begin{align}
  \begin{bmatrix} P-C^T_UC_U & -C^T_U E_1 \\ -E^T_1 C_U & \gamma^2I-E^T_1E_1 \end{bmatrix}
  -\begin{bmatrix}A_U^T  \\ E^T \end{bmatrix}  P \begin{bmatrix} A_U^T  \\ E^T \end{bmatrix}^T>0 .  \label{H_infty_4}
\end{align}
Then, by using Schur complement lemma twice and algebra computaiton, (\ref{H_infty_4}) is equivalent to the following LMI: 
\begin{equation}
    \begin{bmatrix}I \\ A^T \\ B^T \\ E^T \end{bmatrix}^T
    \underbrace{\left[\begin{array}{cc}P^{-1} & 0 \\ 0 & -{G}_{\ast}
    \end{array}\right]  }_{:=G}
    \begin{bmatrix}I \\ A^T \\ B^T \\ E^T \end{bmatrix}    >0    \label{H_infty_5}
\end{equation} 
where
\begin{equation*}{G}_{\ast}=
    {  \begin{bmatrix}I&0\\ K&0 \\ 0&I \end{bmatrix}
    \begin{bmatrix}P-C^T_UC_U& -C^T_U E_1\\ -E^T_1 C_U &\gamma^2 I-E^T_1 E_1\end{bmatrix}^{-1}
    \begin{bmatrix}I&0\\ K&0 \\ 0&I \end{bmatrix} ^T},
\end{equation*}
which is the same condition as (\ref{H_infty_2}).
Then, we can summarize the conclusion:

\textbf{Lemma 2.} Suppose that there exist scalars $\lambda$ and $\beta$, matrices $ Q=Q^T>0$ and $L$ such that (\ref{LMI_H_inf_sum})
holds. Then, the data $(\hat{U},\hat{X},\hat{W})$ can be utilized for $H_{\infty}$ control of closed-loop system with a performance index $\gamma$.

\begin{table*}
    \centering
    \begin{minipage}{1\textwidth}
    \hrule
    \begin{subequations}\label{LMI_H_inf_sum}
    \begin{align}
        & \left[\begin{array}{ccccccc}Q-\beta I & 0 & 0& 0& 0& 0& 0\\
             \ast& 0& 0& 0& Q& 0& 0\\  \ast& \ast& 0& 0& L& 0& 0\\   \ast& \ast& \ast& 0& 0& I& 0\\
             \ast& \ast& \ast& \ast& Q& -C^T_{L}E_1& C^T_{L}\\     \ast& \ast& \ast& \ast& \ast& \gamma^2I-E^T_1E_1& 0\\
             \ast& \ast& \ast& \ast& \ast& \ast& I
        \end{array}\right] +\lambda
        \left[\begin{array}{cc}
            \hat{X}_+\\  -\hat{X}_-\\-\hat{U}\\-\hat{W}\\0\\0\\0
         \end{array} \right]
        \left[\begin{array}{cc}
            \hat{X}_+\\ -\hat{X}_-\\-\hat{U}\\-\hat{W}\\0\\0\\0
         \end{array} \right]^T
        \geq 0   ,                          \label{LMI_H_inf}\\
        &\left[\begin{array}{ccc} Q& -C^T_{L} E_1 & C^T_{L} \\  \ast& \gamma^2I-E^T_1 E_1& 0\\    \ast& \ast& I 
        \end{array}\right]  >0       .     \label{LMI_for_H_inf}
    \end{align}
    \end{subequations}
    \medskip
    \hrule
    \end{minipage}
\end{table*}

\textbf{Proof:}
We first prove the conditions that satisfie the demands of Lemma 1. Set 
\begin{align*}
    G&=\left[\begin{array}{c:c}
        G_{11}& G_{12}  \\ \hdashline G_{12}^T& G_{22}
    \end{array} \right]:=\left[\begin{array}{c:c}
        P^{-1}&0  \\ \hdashline 0 & -{G}_{\ast}
    \end{array} \right]  , \\
    H&=\left[\begin{array}{c|c}
        H_{11}&H_{12} \\\hdashline   H_{12}^T &  H_{22}
    \end{array} \right]
    := \xi  
    \begin{bmatrix}
    0 &0 \\ 0 & -I
    \end{bmatrix}
    \xi ^T  ,
\end{align*} 
where
\begin{align*}
\xi=\left[\begin{array}{c:cccccc}
        I &  0 &  0&   0 \\        \hat{X}_+^T&   -\hat{X}_-^T&   -\hat{U}^T&  -\hat{W}^T
    \end{array} \right]^T .
\end{align*} 
By (\ref{H_infty_4}) and (\ref{H_infty_5}), it is clear that $-{G}_{\ast} \leq 0$, then the condition $G_{22}\leq 0$ of Lemma 1 can be verified.
According to the construction of $H$, we can find that the assumption (2) of Lemma 1 is satisfied.
As to the assumption (3), which can be achieved by utilizing (12) and the construction of $H$. 
Define 
\[F=\begin{bmatrix}
    A^T_{\nu}\\B^T_{\nu}\\E^T_{\nu}
\end{bmatrix}  . \] 
By (\ref{H_infty_5}), $G_{11}+F^{T}G_{22}F$ is equivalent to 
\[P^{-1}-\begin{bmatrix}
    A^T_{\nu}\\B^T_{\nu}\\E^T_{\nu}
\end{bmatrix}^T
{G}_{\ast}
\begin{bmatrix}
    A^T_{\nu}\\B^T_{\nu}\\E^T_{\nu}
\end{bmatrix} > 0 .  \]
Then,
\begin{align*}
    H_{22}F&=\left[\begin{array}{c} -\hat{X}_- \\ -\hat{U}_-\\ -\hat{W}_-  \end{array}\right]
     \left[\begin{array}{c} -\hat{X}_- \\ -\hat{U}_-\\ -\hat{W}_-  \end{array}\right]^T 
     \begin{bmatrix}
        A^T_{\nu}\\B^T_{\nu}\\E^T_{\nu}
    \end{bmatrix},\\
    H_{12}&=\left[\begin{array}{c} -\hat{X}_- \\ -\hat{U}_-\\ -\hat{W}_-  \end{array}\right]-\hat{X}^T_+
\end{align*}
where
\begin{align*}
    -\hat{X}^T_+=\left[\begin{array}{c} -\hat{X}_- \\ -\hat{U}_-\\ -\hat{W}_-  \end{array}\right]^T 
    \begin{bmatrix}
       A^T_{\nu}\\B^T_{\nu}\\E^T_{\nu}
   \end{bmatrix}.
\end{align*}
We can verify $H_{22}F=H^T_{12}$.

Hence, by Lemma 1, one can be concluded that
\begin{align}
    G-\lambda H \geq \begin{bmatrix}\beta I&0\\0&0  \end{bmatrix}  \label{G-lambdaH}
\end{align}
for some $\lambda \in \mathbb{R}$ and $\beta \in \mathbb{R}$.
By Schur complement lemma, (\ref{G-lambdaH}) is equivalent to
\begin{align}
    G^{'}-\lambda H^{'}      \geq   0     \label{G-lambdaH_1}
\end{align}
where 
\begin{align*}
    G^{'}&= \begin{bmatrix}P^{-1}-\beta I & 0 \\ 0 &
        \begin{bmatrix}
        0& \begin{array}{cc}\begin{bmatrix}    
            I&0\\ K&0 \\ 0&I
           \end{bmatrix} &0  \end{array}  \\
           \begin{array}{c}\begin{bmatrix}
            I& K^T & 0\\ 0 & 0 & I
         \end{bmatrix}     \\ 0  \end{array}          & \tilde{G}
        \end{bmatrix}  
    \end{bmatrix}  ,  \\
    H^{'}&=     \xi^{'} 
    \begin{bmatrix}
    0 &0 \\ 0 & -I
    \end{bmatrix}
    \xi^{'T}  ,  \\
\tilde{G} &=
\begin{bmatrix}P & -C^T_U E_1 & C^T_U  \\ -E^T_1 C_U &\gamma^2 I-E^T_1 E_1 & 0 \\C_U &0&I
\end{bmatrix}, \\
\xi^{'} &= \left[\begin{array}{c:cccccc}
    I & 0 & 0 & 0 & 0 & 0 & 0\\ \hat{X}_+^T&   -\hat{X}_-^T&   -\hat{U}^T&  -\hat{W}^T& 0& 0& 0
\end{array} \right]^T.
\end{align*}
By setting $P^{-1}=Q$ and $K=LQ^{-1}$ and
multipling (\ref{G-lambdaH_1}) from the left and right sides by a diagonal matrix $diag \{I, I, I, I, Q, I, I\}$,
the matrix $G^{'}$ can be transformed as follows
\begin{align*}
   \begin{bmatrix}Q-\beta I & 0 & 0& 0& 0& 0& 0\\
    \ast& 0& 0& 0& Q& 0& 0\\  \ast& \ast& 0& 0& L& 0& 0\\   \ast& \ast& \ast& 0& 0& I& 0\\
    \ast& \ast& \ast& \ast& Q& -C^T_{L}E_1& C^T_{L}\\     \ast& \ast& \ast& \ast& \ast& \gamma^2I-E^T_1E_1& 0\\    \ast& \ast& \ast& \ast& \ast& \ast& I
   \end{bmatrix}  
\end{align*}
where $C_{L}=CQ+BL$. This verifies that (\ref{LMI_H_inf_sum}) implies (\ref{H_infty_2}). 
Then, the closed-loop system is internal stable and has $H_{\infty}$ performance level $\gamma$ from $w(t)$ to $y_1(t)$.
$\square$  \\

In next subsection, we will present the Lagrange dual formulation of the minimax optimization and use the above conclusion to the formulation. 

\subsection{Lagrange dual formulation of minimax optimizaiton}

The foundation of MPC involves the real-time solution of a constrained optimization problem, which is updated at each sampling interval based on the current state. 
Before the next sampling instant, the resulting control input is applied to the acutal system. 
In the context of robust MPC, by utilizing the current state $x(t)$ in a moving horizon approach,
{our objective is generally to address a minimax optimal control problem of the system (\ref{sys_1_sum}),}
\begin{equation}\min_{u\in\mathbb{U}}\max_{w\in\mathbb{W}}\sum_{i=t}^{\infty}\|y_{1}(i)\|^{2},     \label{minmax_1} \end{equation}
where $\mathbb{U}$ and $\mathbb{W}$ represent the set of all considerable controls and permissible disturbances, respectively.
In the subsequent analysis, a Lagrange duality is employed to obtain a approximation of the minimax problem (\ref{minmax_1}) that is easier to solve.

According to Assumption 1, the allowable disturbances in system (1) can be presented as
\begin{equation}
    \mathbb{W}=\left\{{w\in\mathbb{R}^{p_1\times [0,\infty)}}\bigg|\sum_{i=0}^\infty\|w(i)\|^2\leq\alpha\right\}. \label{noise_1} 
\end{equation}
By (\ref{noise_1}), we can construct a Lagrangian by combining the original objective function and constraints for any $u\in \mathbb{U}$ and $w\in \mathbb{W}$
as follows
\begin{align}
    \sum_{i=t}^\infty\|y_1(i)\|^2& \leq\sum_{i=t}^{\infty}(\|y_{1}(i)\|^{2}-\gamma^{2}\|w(i)\|^{2})+\gamma^{2}\alpha  \notag \\
    &\leq\max_{w\in L_2}\left(\sum_{i=t}^{\infty}\|y_1(i)\|^2-\gamma^2\|w(i)\|^2+\gamma^2\alpha\right)  \label{Lagran}
\end{align}
for a constant $\gamma> 0$. Then, define a function
\begin{equation}V(x)=\max_{w\in L_2}\left(\sum_{i=0}^\infty(\|y_1(i)\|^2-\gamma^2\|w(i)\|^2)\right)  \label{value_F}   \end{equation}
for the system (\ref{sys_1_sum}), where $x=x(0)$. By using dynamic programming, it can be derived that
\begin{align} &    V(x(t))           \notag     \\
     =&\max_{w(t)\in\mathbb{R}^{p_1}}\left(\|y_1(t)\|^2
    -\gamma^2\|w(t)\|^2+V(x(t+1))\right) .                       
\end{align}
Moreover, we can conclude that system (\ref{sys_1_sum}) satisfies the dissipativity inequality
\begin{equation}V(x (t+1))-V(x(t))\leq\gamma^2\|w(t)\|^2-\|y_1(t)\|^2    \label{diss_ineq}\end{equation}
for $w(t)\in L_2[0,~\infty)$. In the light of (\ref{value_F}), one has from (\ref{Lagran}) that
\begin{equation}
    \sum_{i=t}^\infty\|y_1(i)\|^2\leq V(x(t))+\gamma^2\alpha,
\end{equation}
for any $u\in \mathbb{U}$ and $w\in \mathbb{W}$. Hence, we can arrive at
\begin{align}
    V(x(t))+\gamma^{2}\alpha & \geq\max_{w\in \mathbb{W}}\sum_{i=t}^{\infty}\|y_{1}(i)\|^{2}  \notag \\
    &\geq\min_{u\in \mathbb{U}}\max_{w\in \mathbb{W}}\sum_{i=t}^{\infty}\|y_{1}(i)\|^{2}.   \label{before_opt}
    \end{align}
Then, we can determine the optimal upper bound for (\ref{before_opt}) by the following minimization problem 
\begin{equation}\min_{\gamma^2}V(x(t))+\gamma^2\alpha\quad\text{s. t. (\ref{diss_ineq}) for the system (\ref{sys_1_sum}),} \label{dual_Prob}
\end{equation}
which serves as a Lagrange dual formulation for minimax optimization problem (\ref{minmax_1}).

Let us focus on Lagrange dual formulation (\ref{dual_Prob}). 
It is worth noting that the condition (\ref{H_infty_2}) precisely imply (\ref{diss_ineq}) in the scenario of $u(t)=K x(t)$
by defining $V(t)=x^T(t)Px(t)$ with a positive definite matrix $P$. Hence, we can say that optimization problem (\ref{dual_Prob}) is equivalent to that as follows
\begin{equation}
    \min_{\lambda,\beta,\gamma^2}V(x(t))+\gamma^2\alpha\quad\text{s. t. (\ref{LMI_H_inf_sum}) for the system (\ref{sys_1_sum}),} \label{dual_Prob_1}
\end{equation}

Furthermore, we can derive the following result for $H_{\infty}$ control using (\ref{dual_Prob_1})
by temporarily disregarding the time-domain constraints.

\textbf{Lemma 3.} For given offline data $(\hat{U},~\hat{X},~\hat{W})$ generated by system (1) and external disturbance satisfying Assumption~1,
suppose that the LMI-based optimization problem (\ref{dual_Prob_1}) has an optimal solution ($\lambda_{opt}$,~$\beta_{opt}$,~$\gamma_{opt}$,~$Q_{opt}$,~$L_{opt}$), 
then closed-loop system is internally asymptotically stable under state feedback law given by $K_{opt}=L_{opt}Q^{-1}_{opt}$ 
and achieves an $H_\infty$ performance level of at most $\gamma_{opt}$ from $w$ to $y_1$. 

\textbf{Proof:} The proof follows directly from the discussion above. $\square$ \\

To investigate the scenario with the time-domain constraint (\ref{z_constraint}), we define an ellipsoid for the state $x(t)$ as
\begin{equation}\Psi({P},\sigma_s):=\begin{Bmatrix}x\in\mathbb{R}^n|x^T{P}x\leq \sigma_s\end{Bmatrix}\end{equation}
using a matrix $P$ and a scalar $\sigma_s> 0$.
Then, the formulation for the output constraints, subject to $x(t)\in \Psi ({P},\sigma_s)$, can be presented as 
\begin{align}
    & \quad \max_{t\geq 0}|y_{2m}(t)|^{2}  \notag  \\
    & =    \max_{t\geq 0} x^T(t) C^T_{2,K}  C_{2,K} x(t)      \notag  \\
    &\leq\max_{x\in \Psi} x^T  C^T_{2,K} C_{2,K}x     \notag  \\
    &\leq y_{2m,\max}^{2} , \quad m=1,2,\cdots, q_2     \label{z_max_1}
\end{align}
where $C_{2,K}=C_2+D_2K$. In fact, (\ref{z_max_1}) implies all $x$ satisfying $V(x)=x^T P x\leq \sigma_s$ such that 
$x^T C_{2,K}^T  C_{2,K} x \leq y^2_{2m,~max}$.
By utilizing S-Lemma, we have
\begin{equation}\begin{cases}y_{2i,\max}^{2}-x^{{T}}C_{2,K, i}^{{T}}C_{2,K, i} x
        -\sigma_s \varphi +\varphi x^{{T}}Px\geqslant0,\\i=1,2,\cdots,p_{2},\end{cases} \label{z_max_2}
\end{equation}
for a scalar $\varphi>0$, then (\ref{z_max_1}) holds. For simplicity, suppose $\varphi=\frac{y^2_{2i,max}}{\sigma_s}$,
then the condition that makes formula  (\ref{z_max_2}) holds for all non-zero $x$ is 
$\frac{y_{2i,\max}^{2}}{\sigma_s}P-C_{2,K,i}^{T}C_{2,K,i}\geq 0$ or
\begin{equation}
    \frac{y_{2i, \mathrm{max}}^{2}}{\sigma_s}Q - QC_{2, K, i}^{T}C_{2, K, i}Q \geqslant0. \label{z_max_3}
\end{equation}
With the help of Schur complement lemma, the condition that makes formula (\ref{z_max_3}) holds is that suppose symmetric matrix $\Lambda$ such that
\begin{equation}
    \begin{cases}
    \begin{pmatrix}\frac{1}{\sigma_s}\Lambda&C_{2}Q + D_{2}L\\( C_{2}Q + D_{2}L)^{{T}}&Q
    \end{pmatrix}\geqslant0 ,\\\Lambda_{ii}\leqslant y_{2i,\max}^{2}, i = 1,2,\cdots,q_{2},    \label{LMI_output}
\end{cases}
\end{equation}
holds for $Q$ and $L$.
Similarly, by Schur complement lemma, finding the minimum lower bound for $x(t)^T P x(t)\leq \sigma_s $ corresponds to minimizing $\sigma$ under the constraint
\begin{equation}\begin{pmatrix}\sigma&x(t)^T\\x(t)& {Q}\end{pmatrix}\geq0.  \label{LMI_ellipsoid}
\end{equation}
Therefore, (\ref{dual_Prob_1}) with time-domain constraint becomes
\begin{align}&\min_{\sigma,\lambda,\beta,\gamma^2, {L}, {Q}}\sigma+\alpha \gamma^2 \notag \\
    &\text{s. t. (\ref{LMI_H_inf_sum}), (\ref{LMI_output}), (\ref{LMI_ellipsoid}) and }\sigma\leq \sigma_s. \label{dual_Prob_2}
\end{align}
For the sake of generality, we utilize two weight parameters $r_1$ and $r_2$ and conclude the folowing result based on (\ref{dual_Prob_2}). 

\textbf{Theorem 1:} For a given scalar $\sigma_s>0$ and a matrix $\Lambda$, consider offline data $(\hat{U},~\hat{X}_-,~\hat{W})$ generated by system (1) and external disturbance satisfying Assumption 1. 
If the LMI-based optimization problem 
\begin{align}&\min_{\sigma,\lambda,\beta,\gamma^2, {L}, {Q}}r_1\sigma+r_2\gamma^2  \notag \\
    &\quad\mathrm{s.~t.~(\ref{LMI_H_inf_sum}),~(\ref{LMI_output}),~(\ref{LMI_ellipsoid})}  
    \mathrm{~and~}\sigma\leq \sigma_s\label{dual_Prob_3}
\end{align}
has an optimal solution ($\sigma_t$,~$\lambda_t$,~$\beta_t$,~$\gamma_t$,~$Q_t$,~$L_t$),
then the following properties of the closed-loop system hold:\\
(i) the closed-loop system is internally asymptotically stable under state feedback $K_{opt}=L_{opt}Q^{-1}_{opt}$ \\
(ii) the closed-loop system achieves the optimal $H_\infty$ performance level $\gamma_{opt}$ from $w$ to $y_1$. \\
(iii) the constraint (\ref{z_constraint}) is satisfied for $t\geq 0$.

\textbf{Proof:} On the basis of Lemma 2, we can easily obtain the conclusion (i) and (ii). For conclusion (iii), it can be derive the 
closed-loop system satisfies dissipation inequality (\ref{diss_ineq}) that solution ($\sigma_t$,~$\lambda_t$,~$\beta_t$,~$\gamma_t$,~$Q_t$,~$L_t$)
satisfy LMI (\ref{LMI_for_H_inf}) and $V(x)\geq 0$.     \hfill$\Box$ \\

\textbf{Remark 2:}
It is worth noting that the two weighting parameters $r_1$, $r_2$ introduced in (\ref{dual_Prob_3}) allow more flexibility in adapting this optimization problem to different scenarios.
If one wants a larger weight ratio of $\sigma$, i.e., if one values the optimization effect of $\sigma$ more, one can set the parameter $r_1$ larger, 
and vice versa. In addition, the optimization problem (\ref{dual_Prob_3}) can transform into different optimization problems depending on the value of $r_1$ and $r_2$,
e.g., optimization problem (\ref{dual_Prob_3}) can be transformed into the problem (\ref{dual_Prob_2}) by setting $r_1=1$ and $r_2=\alpha$, 
optimization problem (\ref{dual_Prob_3}) implies the cost function in \cite{kothare1996robust} by setting $r_1=1$ and $r_2=0$, 
and optimization problem (\ref{dual_Prob_3}) corresponds to the cost function in \cite{CHEN20061033} by setting $r_1=0$ and $r_2=1$.

\textbf{Remark 3:}
Note that the inequality (10) implies the dissipation inequality (22), 
indicating that (22) can be derived not only as shown in this paper but also from (10) (i.e., (8) and (9)). 
This method is detailed in \cite{CHEN20061033}. Therefore, this paper leverages their correlation to propose a data-driven \(H_{\infty}\) optimal control method.

\textbf{Remark 4:}
In minmax MPC, feedback predictive control is typically employed to mitigate the effects of disturbances and uncertainties 
while also reducing computational complexity \cite{xie2024data}. However, implementing feedback minmax MPC may encounter feasibility and practical issues \cite{chen2006improved}. 
To address this, this paper leverages Lagrange duality to relax the max-min problem into a more tractable minimization problem. 
This approach not only alleviates infeasibility to some extent but also reduces the conservatism introduced by ellipsoid evaluation 
in time-domain constraints.

\subsection{Moving horizon predictive control}
For the moving horizon predictive control, the optimization problem (\ref{dual_Prob_3}) using current state $x(t)$ will be solved in real-time for every time instants $t\geq 0$.
This implementation allows the current state variable $x(t)$ to be used for achieving feedback control. 
However, the dissipation of closed-loop systems may not be guaranteed under the moving horizon optimization control. 
To address this problem, in light of \cite{Schere1185016}, we can add the dissipative constraint condition to the optimization problem (\ref{dual_Prob_3}) as 
\begin{equation}\begin{pmatrix}{x}(t)^T {P}_{t-1} {x}(t) +p_0-p_{t-1} &  {x}(t)^T\\
     {x}(t)& {Q}\end{pmatrix}\geq0 \label{dissi_constraint}
\end{equation}
for each time instants $t\geq 0$, where $p_0=x(0)^T P_0 x(0)$ and $p_t$ is calculated from equation:
\begin{equation}p_t:=p_{t-1}-\big[{x}(t)^T{P}_{t-1} {x}(t)- {x}(t)^T {P}_t {x}(t)\big]. \label{p_t}
\end{equation}                                    
Then the optimization problem (\ref{dual_Prob_3}) is formulated as 
\begin{align}&\quad\quad\quad\min_{\sigma,\lambda,\beta,\gamma^2, {L}, {Q}} r_1 \sigma +r_2\gamma^2   \notag \\
    &\quad\text{s.~t.~(\ref{LMI_H_inf_sum}),~(\ref{LMI_output}),~(\ref{LMI_ellipsoid}),~(\ref{dissi_constraint})} 
    \text{~and~} 
    \sigma \leq \sigma_s  \label{dual_Prob_4}
\end{align}
for each time $t>0$.

By \cite{CHEN20061033}, we know that a prerequisite for ensuring that the closed-loop system satisfies constraint (\ref{z_constraint}) is that 
the initial state of the system belongs to the elliptical domain. 
Therefore, in case of unsolvable situations, we can avoid this problem by continuously increasing the value of $\sigma_s$.
By incorporating a scalar $\eta \geq 0$, we can reformulate the optimization problem (\ref{dual_Prob_4}) into a more feasible version
\begin{align}&~~~\min_{\sigma,\lambda,\beta,\gamma^2, {L}, {Q}} r_{1} \sigma  +r_{2}\gamma^{2}\notag \\  
    &\text{s.~t. (\ref{LMI_H_inf_sum}), (\ref{LMI_output}), 
    (\ref{LMI_ellipsoid}), (\ref{dissi_constraint})},~
    \mathrm{and}~\sigma  \leq \sigma_{s}(1+\eta).  \label{dual_Prob_5}
\end{align}
Furthermore, the following feasibility conclusion can be reached from the above analyses. \\
\textbf{Lemma 4:} For given $\sigma_s> 0$ and $\Lambda$, suppose that the Assumption 1 holds and 
LMIs (\ref{LMI_H_inf}), (\ref{LMI_for_H_inf}), (\ref{LMI_ellipsoid}) and (\ref{LMI_output}) are feasible with a solution ($\sigma$,~$\lambda$,~$\beta$,~$\gamma$,~$Q$,~$L$) 
at time $t_0\geq 0$. Then, the feasibility of the optimization problem (\ref{dual_Prob_5}) can be guaranteed at every time $t_0+n$ for some $n\geq 0$ and $\eta \geq 0$.\\

\noindent \textbf{Proof}:  
For $t=0$, the feasibility of the optimization problem (\ref{dual_Prob_3}) implies the feasibility of the optimization problem (\ref{dual_Prob_4}).
Let there exists a bounded initial state $x(0)$ such that (\ref{LMI_ellipsoid}) holds by defining $\sigma_0=x^T(0)P_0x(0)$.
On the basis of the fact that (\ref{LMI_H_inf_sum}) and (\ref{LMI_output}) do not rely on the system state variable $x(k)$,
it can be concluded that the problem is initially feasible and remains feasible in the future. 
Therefore, (\ref{dual_Prob_3}) is feasible for $\sigma_0\leq \sigma_s$. In addition, the fact (\ref{LMI_H_inf_sum}) is feasible means that
(\ref{diss_ineq}) is feasible with $\gamma=\gamma_0$ and $V(x(t))=x(t)^TP_0x(t)$, and then $x(t)$ is bounded for $t=1$.

When $t>0$, suppose that there exists a bounded system state $x(t)$ such that the optimization problem (\ref{dual_Prob_4}) has a set of feasible solution ($\sigma_t$,~$\lambda_t$,~$\beta_t$,~$\gamma_t$,~$Q_t$,~$L_t$).
Hence, $x(t+1)$ is a bounded system state by using (\ref{diss_ineq}).

For $t=t+1$, let there exists a bounded $\sigma_{t+1}=x(t+1)P_{t}x(t+1)$ such that (\ref{LMI_ellipsoid}) holds. 
Then, the fact (\ref{dissi_constraint}) is feasible at time instant $t$ implies that $p_0-p_{t}\geq 0$ by (\ref{p_t}). 
That is to say, when $t=t+1$, (\ref{dissi_constraint}) is feasible with $Q=Q_t$.
Therefore, the optimization problem (\ref{dual_Prob_4}) has a fesible solution ($\sigma_{t+1}$,~$\lambda_t$,~$\beta_t$,~$\gamma_t$,~$Q_t$,~$L_t$) to the system state variable $x(t+1)$ if $\sigma_{t+1}\leq \sigma_s$.
\hfill$\Box$ 

Herein, the corresponding moving horizon predictive control algorithm is presented as Algorithm 1.

\begin{algorithm}
    \caption{Moving Horizon Control Algorithm}
    \label{alg:1}
    \noindent\textbf{1):} Initialization. Set $t=0$ and given $\sigma_s$, $\Lambda$, initial state $x(0)$ and offline data $(\hat{U},\hat{X},\hat{W})$. \vspace{\baselineskip}

    \noindent \textbf{2):} For $t=0$, find a set of solutions $(\sigma _0, \lambda_0,\beta_0,$$\gamma_0, Q_0, L_0)$ to the optimization problem (\ref{dual_Prob_3}).
    If no feasible solution is found, substitute $\sigma_s$ by $\sigma_s(1+\eta)$ and increase $\eta\geq 0$ slightly each time. 
    Let $t_0=L_0 Q_0, P_0=Q^{-1}_0, p_0=x^T(0) P_0 x(0)$ and go to \textbf{4)}.\vspace{\baselineskip}
    
    \noindent \textbf{3):} When $t > 0$, given $x(t)$ and find a set of solutions $(\sigma_t,\lambda_t,\beta_t, \gamma_t, Q_t, L_t)$ 
    to the optimization problem (\ref{dual_Prob_4}), if not feasible, substitute
    the optimization problem (\ref{dual_Prob_4}) by (\ref{dual_Prob_5}). Set $K_t=L_t Q_t, P_t=Q^{-1}_t$, 
    and $p_t:=p_{t-1}-{x}^T(t){P}_{t-1}x(t)+{x}^T(t){P}_t{x}(t)$.\vspace{\baselineskip}
    
    \noindent \textbf{4):} Achieve control input
    \begin{equation}
    u(t)=K_t x(t),~\forall t\geq 0 \label{controller}
    \end{equation}
    and apply it to the system. Let $t=t+1$, and then proceed to \textbf{3)} continuously.\vspace{\baselineskip}
\end{algorithm}

Then, we now elaborate the following conclusion of the properties of the closed-loop system based on the above discussion.
The conclusion integrates the results of Theorem 1 into the MPC framework, allowing the $ H_{\infty} $ performance to be continuously optimized.
\vspace{\baselineskip}

\noindent \textbf{Theorem 2:} For a given $\sigma_s> 0$, suppose that\vspace{\baselineskip}

$\bullet $ Assumption 1 holds; 

$\bullet $ The offline data set $(\hat{U},~\hat{X},~\hat{W})$ generated by system (1) is admissible;

$\bullet $ The optimization problem (\ref{dual_Prob_5}) is feasible with a solutions ($\lambda_t$,~$\beta_t$,~$\sigma_t$,~$\gamma_t$,~$Q_t$,~$L_t$) at time $t\geq 0$;

$\bullet $ The set of optimal solutions ($\lambda_t$,~$\beta_t$,~$\sigma_t$,~$\gamma_t$,~$Q_t$,~$L_t$) is bounded,

\vspace{\baselineskip}

then the closed-loop system with controller (\ref{controller}) can reach the properties as follows:\vspace{\baselineskip}

(i) The time-domain constraint is satisfied at every time $t\geq 0$;

(ii) The stabilization of the system can be achieved for finite energy disturbances;

(iii) For the discrete-time LTI system (\ref{sys_1_sum}), the $H_\infty$ norm is bounded above by 
$\gamma_\infty:= \lim_{t\rightarrow \infty} \textrm{max}\{\gamma_t\} < \infty$. 

(iv) The closed-loop system exhibits disturbance attenuation, i.e.,
\begin{equation}\sum_{i=t_0}^t\|y_1(i)\|^2-\gamma^2\|w(i)\|^2\leq x(t_0)^T{P}_{t_0}x(t_0)+p_0-p_{t_0}\end{equation}
holds with $\gamma\leq \infty$ and $p_0-p_{t_0}\geq 0$ for any $t\geq t_0\geq 0$. \vspace{\baselineskip}

\noindent \textbf{Proof:} For a given system state $x(t) $, property (i) can be clearly derived from the solvability of (\ref{LMI_output}),
i.e., the system state $x(t) $ and the calculated feedback gain $K$ satisfy the given constraint condition.

Suppose there exist an optimal solution $(\sigma_t,\lambda_t,\beta_t,\gamma_t,$ $ Q_t, L_t)$, 
for any $t\geq t_0 \geq 0$, when $\gamma=\gamma_t$, $V(x)=x^T(t) P_t x(t)$ and $P_t=Q^{-1}_t$, 
(\ref{LMI_H_inf_sum}) can derive (\ref{diss_ineq}) and one has
\begin{align*}
    &\|y_{1}(t_{0})\|^{2}-\gamma_{t_{0}}^{2}\|w(t_{0})\|^{2}    \\
    \leq ~&  {x}(t_0)^T {P}_{t_0} {x}(t_0)- {x}(t_0+1)^T {P}_{t_0} {x}(t_0+1) \\
    &\|y_{1}(t_{0}+1)\|^{2}-\gamma_{t_{0}+1}^{2}\|w(t_{0}+1)\|^{2}  \\
    \leq ~& x(t_0+1)^T {P}_{t_0+1} {x}(t_0+1) \\
    &-x(t_0+2)^T {P}_{t_0+1} {x}(t_0+2)\\
    &~~~~~\cdots\\  
   &\|y_{1}(t)\|^{2}-\gamma_{t}^{2}\|w(t)\|^{2} \\
   \leq ~ & x(t)^{T} {P}_{t}x(t)-x(t+1)^{T} {P}_{t}x(t+1)
\end{align*}
and then
\begin{align}&\sum_{i=t_{0}}^{t} \Big( \|y_{1}(i)\|^{2}-\gamma_{i}^{2}\|w(i)\|^{2} \Big) \notag \\
    &\leq x(t_{0})^{T} {P}_{t_{0}}x(t_{0})-x(t+1)^T {P}_t x(t+1)   \notag \\
    &-\sum_{i=t_{0}+1}^{t} \Big( x(i)^{T} {P}_{i-1}x(i)-x(i)^{T} {P}_{i}x(i)\Big)
     .  \label{index_27}
\end{align}
Furthermore, according to Schur complement lemma, we can yield from the feasibility of (\ref{dissi_constraint}) that  
\begin{equation}p_0-p_{t-1}+x(t)^T{P}_{t-1}x(t)-x(t)^T{P}_t x(t)\geq0. \label{th2_proof_1}  \end{equation}
Substituting (\ref{p_t}) recursively into the inequality (\ref{th2_proof_1}) from time instant $t=0$ to time instant $ t=t_0$,
we can draw from the dissipation constraint condition that
\begin{equation}p_0-p_{t_0}+\sum_{i=t_0+1}^t\left[x(i)^T{P}_{i-1}x(i)-x(i)^T{P}_ix(i)\right]\geq0.\end{equation}
Hence, the inequality (\ref{index_27}) can be reformed as 
\begin{align}\sum_{i=t_0}^t\|y_1(i)\|^2-\gamma^2\|w(i)\|^2&\leq x(t_0)^T{P}_{t_0} {x}(t_0)+p_0-p_{t_0}  \notag \\
    &-x(t+1)^T{P}_t {x}(t+1)  \label{index_30}
\end{align} 
for $\gamma_{max} < \infty$. Because of the positive definiteness of $P_t$, by (\ref{index_30}), we can arrive at conclusion (iv). 
Furthermore, for time $t=t_0$, the solvability of (\ref{dissi_constraint}) can indicate $p_0-p_{t_0}\geq 0$ based on (\ref{p_t}). Then, for the case of $t\rightarrow \infty$,
we can conclude the disturbance satisfying Assumption 1 has the limit and 
\begin{equation}\sum_{i=0}^{\infty}\|y_{1}(i)\|^{2}\leq x(0)^{T}{P}_{0}
    {x}(0)+\gamma_{\infty}^{2}\sum_{i=0}^{\infty}\|w(i)\|^{2}.
\end{equation}    
This implies conclusion (ii). As to conclusion (iii), when $x(0)=0$, we find that
\begin{equation}
    \gamma^{2}\sum_{i=0}^{\infty}\|w(i)\|^{2}\geq \sum_{i=0}^{\infty}\|y_{1}(i)\|^{2} .
\end{equation}
for finited energy disturbance $w(t)$.   \quad \hfill $\Box$

\section{Numerical Example}

In this subsection, we consider a batch reactor system \cite{walsh2001scheduling,de2019formulas} as simulation application to demonstrate the proposed scheme. 

The open-loop unstable system represented by equation (\ref{sys_1_sum}) is characterized by the matrices 
\begin{align*}
A_{\nu }&=\begin{bmatrix}
1.178 & 0.001& 0.511 &  -0.403\\
-0.051& 0.661&-0.011 &  0.061\\
0.076 & 0.335&  0.560&  0.382     \\
0     & 0.335& 0.089 &  0.849
\end{bmatrix},  \\
B_{\nu }&=\begin{bmatrix}   
0.004 & -0.087 \\
0.467 & 0.001\\
0.213 &-0.235\\
0.213 &-0.016
\end{bmatrix}.
\end{align*}
The other system parameters are considered as   
\begin{align*} 
C_1 &=\begin{bmatrix}
    1 &0 &1 &-1
\end{bmatrix},~D_1=0,\\
C_2&=\begin{bmatrix}
    0.5 &0.5 &1& 1   
\end{bmatrix},~D_2=\begin{bmatrix}
    0 &1 
\end{bmatrix}
\end{align*}
and $E_{\nu }=E_1=I$. These specific values indicate the system's inherent instability under open-loop conditions, 
necessitating further analysis and potentially the design of a suitable control strategy to achieve desired performance and stability.

\begin{figure}[ht!]
    \begin{center}
    \includegraphics[height=5.2cm,width=0.5\textwidth]{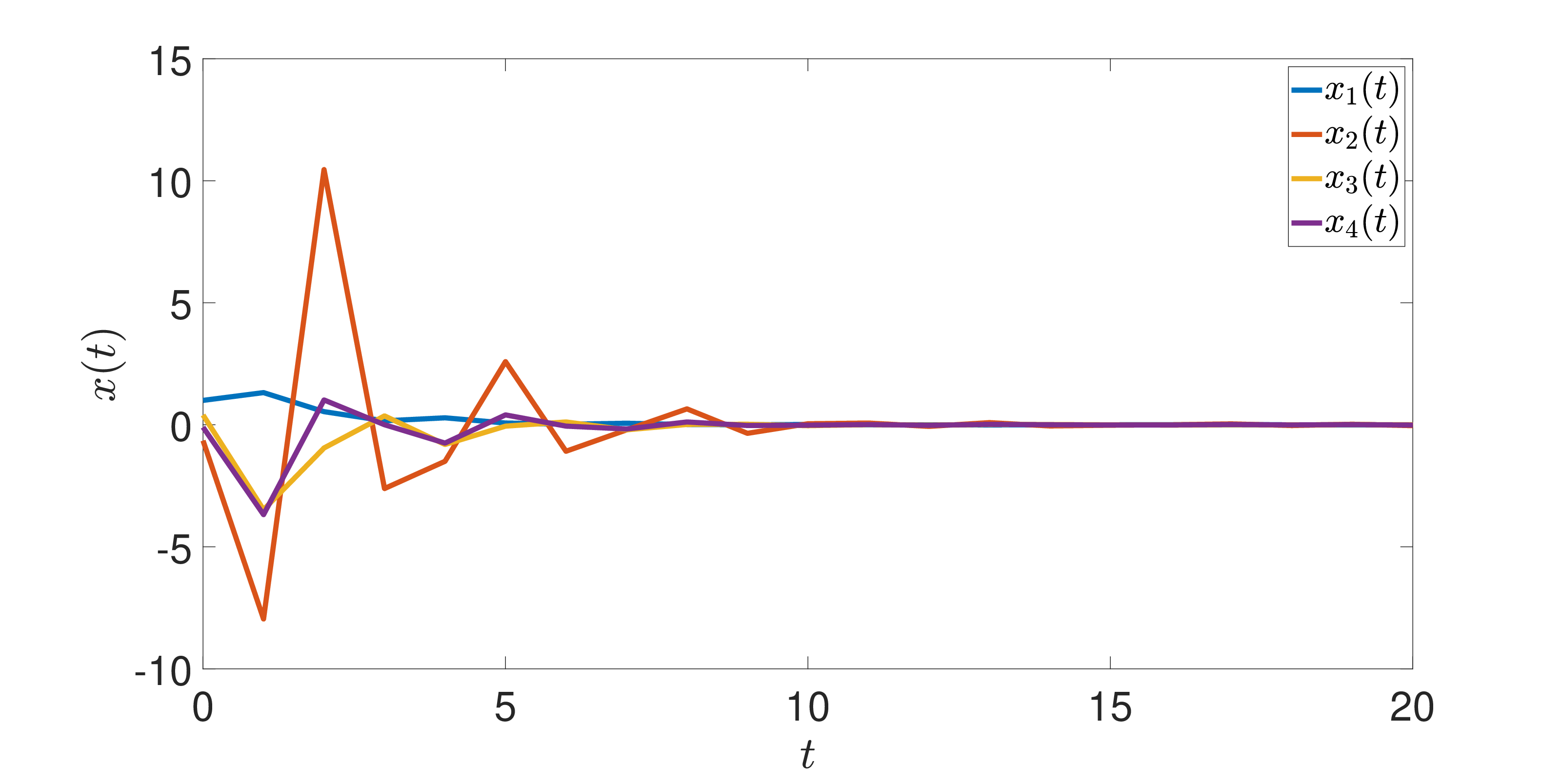}    
    \caption{State trajectories under control input}  
    \label{fig1}                                 
    \end{center}                                 
    \end{figure}

In this example, our objective is to determine a control gain with optimal $H_\infty$ performance and state ellipsoid of an unknown system by utilizing input-state-disturbance data.
The external disturbance $w(t)$ is assumed to satisfy $||w(t)||^2\leq \alpha_{t}$ for all $t$.
Then, we obtain input data over a time length of \( T = 20 \) and initial state data from a standard normal distribution;
the external disturbance data are randomly sampled to adhere to the earlier assumption and standard normal distribution.  
Simultaneously, the corresponding state data can be generated. 
Herein, we implement the proposed data-driven control scheme to system (\ref{sys_1_sum}) based on the off-line data $(\hat{U},~\hat{X},~\hat{W})$ and parameters $\alpha_{t}=10^{-4}$,
$\sigma_s=10$, $\Lambda=1.2I$, $r_1=1$, $r_2=1$, initial condition $x(0)=[1;~-0.65;~0.4;~-0.1]$ and random external disturbance $w(t)$. 
The optimization problem is solved by Yalmip interface with Mosek solver in MATLAB.
We can find a stabilizing controller for system (1) by using the (\ref{controller}) as displayed in the Figure \ref{fig1}. From the figure, the fact is clear that
the state responses of closed-loop system converge to equilibrium points under effectiveness of the control input.

    \begin{figure}[ht!]
        \begin{center}
        \includegraphics[height=5.2cm,width=0.5\textwidth]{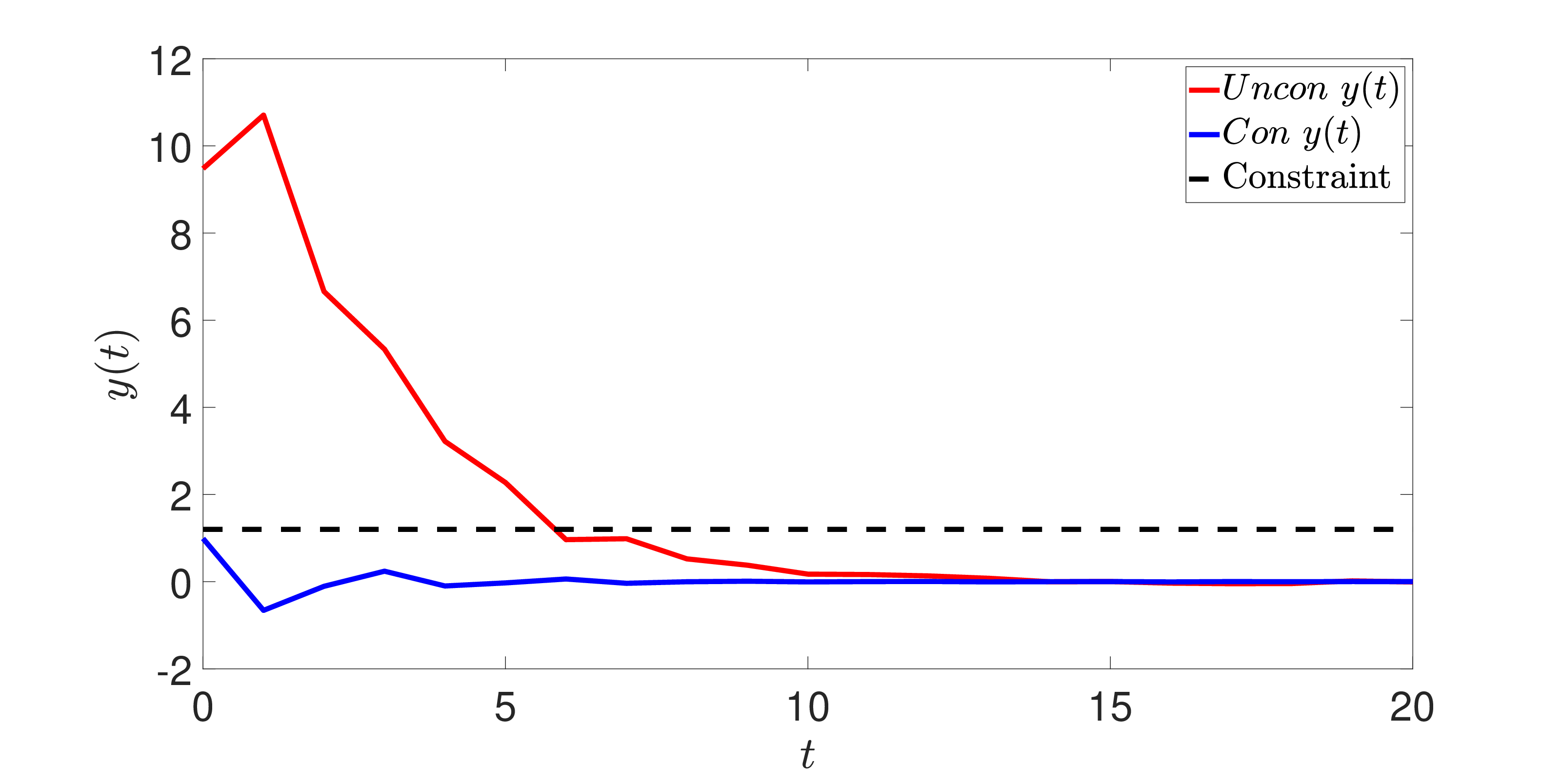}    
        \caption{The comparison of control output between the unconstrained control and the constrained method}  
        \label{fig3}                                 
        \end{center}                                 
    \end{figure}
\begin{figure}[ht!]
    \begin{center}
    \includegraphics[height=5.2cm,width=0.5\textwidth]{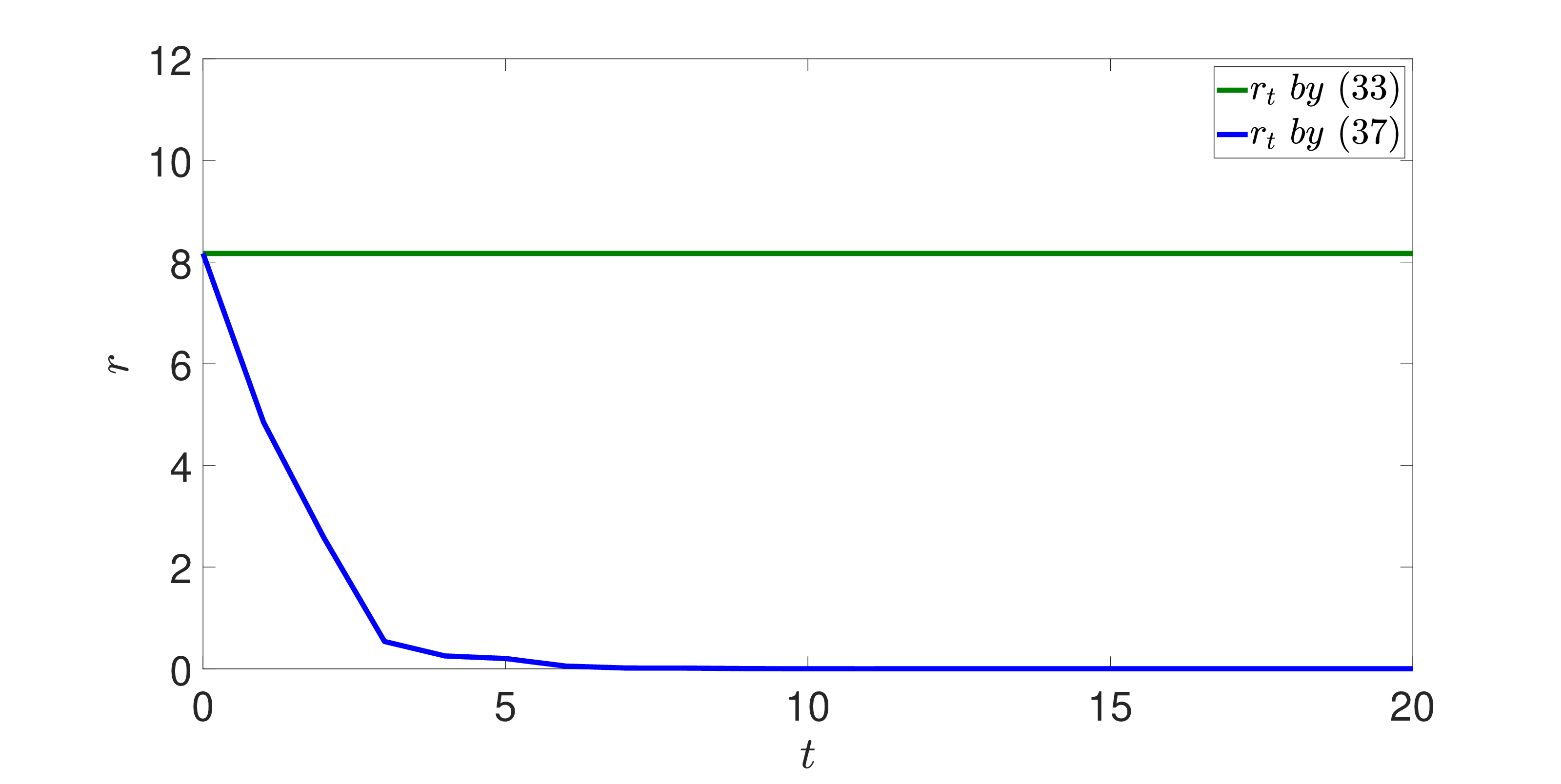}    
    \caption{The comparison of $r_t$ between moving horizion method and static method   }  
    \label{fig4}                                 
    \end{center}                                 
\end{figure}

\begin{figure}[ht!]
    \begin{center}
    \includegraphics[height=5.2cm,width=0.5\textwidth]{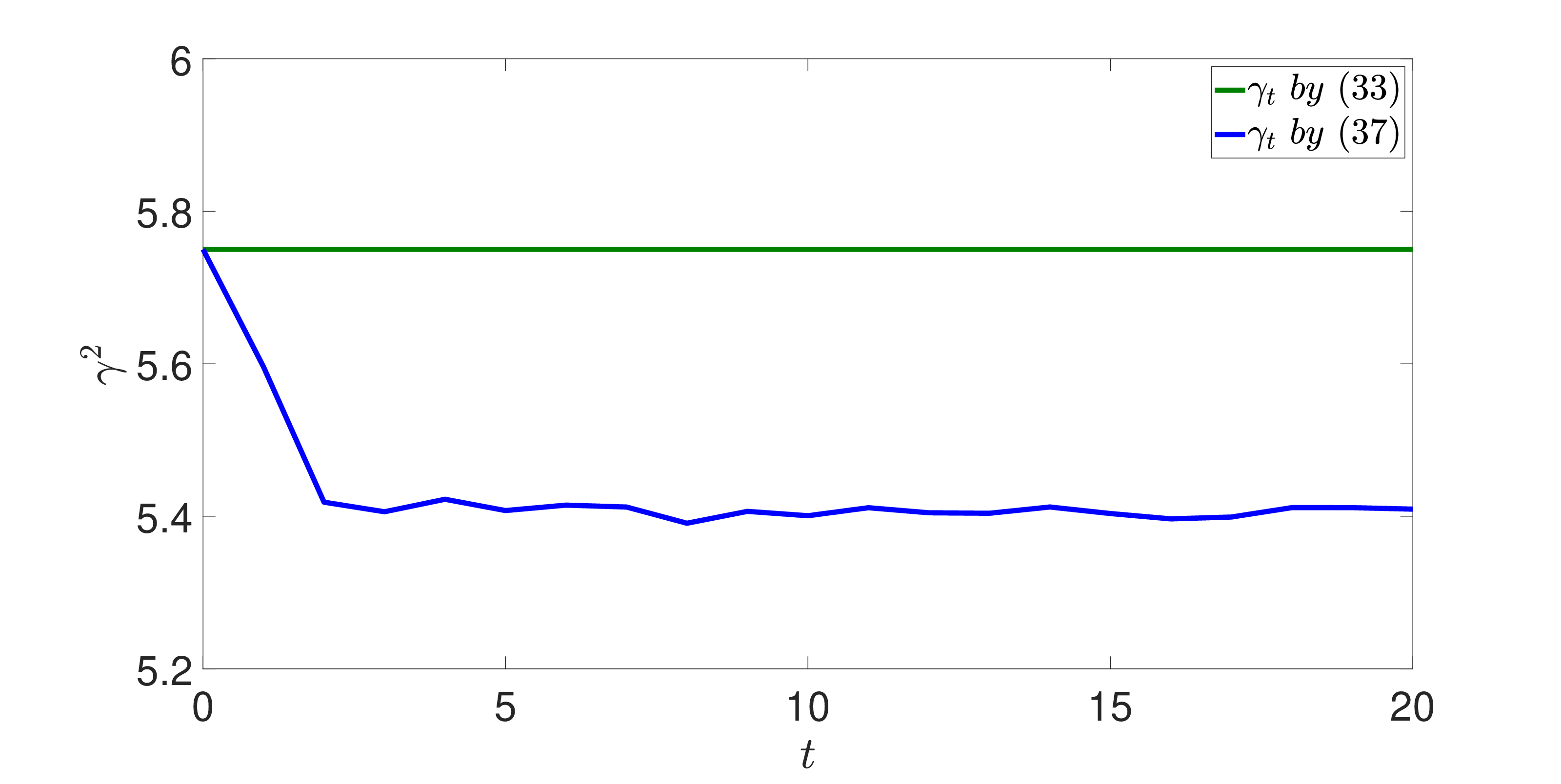}    
    \caption{The comparison of $\gamma_t$ between moving horizion method and static method}  
    \label{fig5}                                 
    \end{center}                                 
\end{figure}
In Figure \ref{fig3}, we can find that the control output curve ($Uncon\ z(t)$) without the LMI constraint (\ref{LMI_output}) exceeds the allowable limit,
while the control output curve ($Con\ z(t)$) using the proposed method stays within the constraint. This demonstrates that the implementation of 
moving horizon predictive control (\ref{dual_Prob_5}) does not violate the LMI constraints. 
Then, we illustrate the superiority of the optimization problem (\ref{dual_Prob_5}) with moving horizon control, compared to the problem (\ref{dual_Prob_3}) 
without horizon control as shown in Figure~\ref{fig4} and Figure~\ref{fig5}. 
More specifically, both methods achieve the same performance levels initially, but as time progresses, the performance levels of ($r_t$ by (\ref{dual_Prob_5})) 
and ($\gamma_t$ by (\ref{dual_Prob_5})) remain consistently optimal. 
Therefore, the moving horizon optimization problem (\ref{dual_Prob_5}) effectively improves performance levels while ensuring compliance with constraints.
In addition, a smaller $\gamma_t$ indicates better disturbance attenuation.

For comparison, the data-driven control schemes from \cite{van2020noisy} and \cite{li2024data} are substituted into the LMI (\ref{LMI_H_inf_sum}) of optimization problem (\ref{dual_Prob_5}). 
However, 
when applying the $H_{\infty}$ design methods from \cite{van2020noisy} and \cite{li2024data} to the optimization problem (\ref{dual_Prob_5})  
using the preceding system matrices and parameters, no feasible solution $(\sigma,\lambda,\beta,\gamma,Q,L)$ was found.
This indicates that the proposed method is less conservative than the methods of \cite{van2020noisy,li2024data}.

There are two reasons for the low conservativeness: First, the output $y(t)$ in this paper, which comes from the state equation with 
$w(t)$, differs from the $y(t)$ in \cite{van2020noisy,li2024data}. In \cite{van2020noisy,li2024data}, $y(t)=Cx(t)+Du(t)$, while in this paper, $y(t)=Cx(t)+Du(t)+Ew(t)$. 
The second reason is that the external disturbance data could have been directly (assumed) utilized, whereas \cite{van2020noisy,li2024data} only requires 
$\|w(t)\|\leq \alpha$, making the system addressed in the proposed method more precise. 
However, we note that the assumption of being able to use external disturbance data is strong while realizable. 
Specifically, one can design a disturbance estimator to estimate the disturbance values of trajectory data, 
for instance, using the approach presented in Section IV of \cite{pan2022stochastic}.
Additionally, it is worth noting that the system models treated in this paper and \cite{van2020noisy,li2024data} are different. 
The state-space equation in this paper includes an additional $E$ matrix in the $w(t)$ term compared to the equation in \cite{van2020noisy,li2024data}, 
which allows the external disturbance to enter the state $x(t)$ from any direction. 
For example, if $E$ is set to the identity matrix $I$, the external disturbance 
$w(t)$ will affect the state $x(t)$ from all directions. Conversely, when $E$ is $[0;0;1]$, the external disturbance 
$w(t)$ will only influence the state component $x_3(t)$ (the choice of the dimensions of the matrix $E$ 
here is for explanatory purposes only). 
Therefore, restricting the direction where the disturbance can enter the system reduces the conservativeness of the control method.

\section{Conclusion}
This article has studied a data-driven $H_{\infty}$ predictive control scheme designed for an unknown system subject to time-domain constraints. 
By leveraging Lagrange duality, the approach transformed the minimax problem into a more tractable minimization problem, 
thus reducing the conservatism associated with ellipsoidal evaluations of time-domain constraints. Utilizing both input-output data and noisy data, 
our scheme has achieved $H_{\infty}$ performance in the closed-loop system. 
The comprehensive analysis conducted demonstrates that the proposed control method has ensured closed-loop stability, effective disturbance attenuation, 
and satisfaction of constraints. 
The validity and advantage of the approach have been further confirmed through numerical simulations involving a batch reactor system, 
highlighting its robustness and feasibility in noisy environments.
This work contributes to the field of data-driven control by offering a robust and practical control strategy that can be applied to various real-world systems.
Future research can build upon this framework by exploring its application to more complex systems and enhancing its adaptability and performance through advanced data-driven techniques.

\bibliographystyle{IEEEtran} 
\bibliography{ref.bib}

\vfill

\end{document}